\documentclass{amsart}

\DeclareMathSymbol{\twoheadrightarrow}
{\mathrel}{AMSa}{"10}

\def\Q{{\mathbf Q}}
\def\Z{{\mathbf Z}}
\def\C{{\mathbf C}}
\def\R{{\mathbf R}}
\def\F{{\mathbf F}}
\def\L{{\mathbf L}}

\def\SS{{\mathbf S}}
\def\RR{{\mathfrak R}}
\def\A{{\mathbf A}}
\def\Sn{{\mathbf S}_n}
\def\An{{\mathbf A}_n}

    \def\red{\mathrm{red}}
\def\Gal{\mathrm{Gal}}
\def\Perm{\mathrm{Perm}}

\def\Sz{\mathrm{Sz}}
\def\U{\mathrm{U}}

\def\End{\mathrm{End}}
\def\Aut{\mathrm{Aut}}

\def\I{\mathrm{Id}}

\def\fchar{\mathrm{char}}

\def\GL{\mathrm{GL}}
\def\PGL{\mathrm{PGL}}
\def\PSL{\mathrm{PSL}}
\def\SL{\mathrm{SL}}

\def\M{\mathrm{M}}

\def\dim{\mathrm{dim}}

\def\P{{\mathbf P}}

\def\T{{T}}

\newtheorem{thm}{Theorem}[section]
\newtheorem{lem}[thm]{Lemma}

\theoremstyle{definition}
\newtheorem{defn}[thm]{Definition}
\newtheorem{ex}[thm]{Example}
\newtheorem{exs}[thm]{Examples}

\newtheorem{rem}[thm]{Remark}

\title[Hyperelliptic jacobians]
{Hyperelliptic jacobians without complex multiplication, doubly
transitive permutation groups and projective representations}
\author[Yuri G. Zarhin]{Yuri G. Zarhin}
\thanks{Partially supported by the NSF}
\begin{document}
\maketitle

\section{Introduction}

  Let $K$ be a  field of characteristic zero, $K_a$ its algebraic closure, $n\ge 5$ an integer, $f(x)$ an irreducible
polynomial over $K$ of degree $n$, whose Galois group $\Gal(f)$
acts doubly transitively on the set $\RR$ of roots of $f$.  Let
$C:y^2=f(x)$ be the corresponding hyperelliptic curve and X=$J(C)$
its jacobian defined over $K$. Earlier, the author
\cite{ZarhinMRL}, \cite{ZarhinTexel}, \cite{ZarhinMMJ} has proven
that the ring $\End(X)$ of all $K_a$-endomorphisms coincides with
$\Z$ if  $\Gal(f)$ is either the full symmetric group $\Sn$ or the
corresponding alternating group $\An$ or  a {\sl small} Mathieu
group $\M_n$ (with $n=11$ or $12$) or $\RR$ could be identified
with the projective space $\P^{m-1}(\F_q)$ over a finite field
$\F_q$ of {\sl odd} characteristic in such a way that $\Gal(f)$
contains the projective special linear group $\PSL_m(\F_q)$ while
$m \ge 3$ and $(m,q) \ne (4,3)$. (Similar results were obtained
when $\Gal(f)=\L_2(2^r), \Sz(2^{2r+1})$ or $\U_3(2^r)$
\cite{ZarhinTexel}, \cite{ZarhinPAMS}.) The proof was based on an
observation that in all these cases the natural (faithful)
representation of $\Gal(f)$ in $X_2$ is {\sl very simple}; in
particular, it is {\sl absolutely irreducible}. (See
\cite{ZarhinTexel}, \cite{ZarhinMRL2}, \cite{ZarhinCrelle}
 for the definition and basic properties of very
simple representations.)

We refer the reader to \cite{Mori1}, \cite{Mori2}, \cite{Katz1},
\cite{Katz2}, \cite{Masser}, \cite{ZarhinMRL}, \cite{ZarhinTexel},
\cite{ZarhinMRL2}, \cite{ZarhinPAMS}, \cite{ZarhinMMJ} for a
discussion of known results about, and examples of, hyperelliptic
jacobians without complex multiplication.

In the present paper we suggest a new approach to already known
examples  when $\End(X)=\Z$ from \cite{ZarhinMRL},
\cite{ZarhinTexel} and \cite{ZarhinMMJ}. Namely, instead of very
simplicity, we use a theorem of Feit-Tits \cite{FT} complemented
by results of Kleidman-Liebeck \cite{KL}. Besides obtaining new
proofs of already known results, we get new examples when
$\End(X)=\Z$ with {\sl reducible} (but still absolutely
indecomposable) $\Gal(f)$-module $X_2$. Namely, we prove that
$\End(X)=\Z$ when $\Gal(f)$ is a {\sl big} Mathieu group (with
$n=22$, $23$ or $24$) or $\RR$ could be identified with the
projective space $\P^{m-1}(\F_{q})$ over a finite field $\F_{q}$
of characteristic $2$ in such a way that $\Gal(f)$ becomes either
the projective special linear group $\L_m(q):=\PSL_m(\F_q)$ or the
projective linear group $\PGL_m(\F_q)$ with $m>2$
(except $(m,q)=(4,2), (3,4)$).
We refer to Theorem \ref{jac} (and Definitions \ref{ubd}  and
\ref{cvr}) for a justification for the long title of the present
article.

The paper is organized as follows. In \S \ref{mainr} we state the main results. Section \ref{gt} contains
 auxiliary results from representation theory.
  In \S \ref{AV}  we discuss fields of definitions for endomorphisms of abelian varieties.
  Section \ref{HJ} contains proofs of the main results.

\section{Main results}
\label{mainr} Throughout this paper we assume that $K$ is a field
of characteristic $0$. We fix its algebraic closure $K_a$ and
write $\Gal(K)$ for the absolute Galois group $\Aut(K_a/K)$.

\begin{thm}
\label{mainM} Let $K$ be a field of characteristic zero,
 $K_a$ its algebraic closure,
$f(x) \in K[x]$ an irreducible polynomial of degree $n$. Suppose
$n=11, 12, 22, 23$ or $24$ and the Galois group of $f$ is the
corresponding Mathieu group $\M_n$.
 Let $C_f$ be the hyperelliptic curve $y^2=f(x)$. Let $J(C_f)$ be
its jacobian, $\End(J(C_f))$ the ring of $K_a$-endomorphisms of
$J(C_f)$. Then  $\End(J(C_f))=\Z$.
\end{thm}

\begin{rem}
The case of small Mathieu groups was done in Th. 7.13 on p. 489 of
\cite{ZarhinTexel} (see also \cite{ZarhinMMJ}). However, in this
paper we give a unified proof for all  Mathieu groups.
\end{rem}

\begin{thm}
\label{mainL} Let $K$ be a field of characteristic zero,
 $K_a$ its algebraic closure,
$f(x) \in K[x]$ an irreducible polynomial of degree $n$ and
$\RR\subset K_a$ the set of its roots. We write $\Gal(f)$ for the
Galois group of $f$. Suppose there exist integers $m>2$ and $r \ge
1$ such that $n=\frac{q^m-1}{q-1}$ where $q=2^r$ and $(m,q)\ne
(2,2), (4,2),(3,4)$. Assume, in addition, that $\RR$ could be
identified with the projective space $\P^{m-1}(\F_{q})$ over the
finite field $\F_{q}$ in such a way that $\Gal(f)$ contains
 $\L_m(q):=\PSL_m(\F_{q})$ as a subgroup.
(E.g., $\Gal(f)=\L_m(q)$ or $\PGL_m(\F_{q})$.)
 Let $C_f$ be the
hyperelliptic curve $y^2=f(x)$. Let $J(C_f)$ be its jacobian,
$\End(J(C_f))$ the ring of $K_a$-endomorphisms of $J(C_f)$. Then
$\End(J(C_f))=\Z$.
\end{thm}

\begin{rem}
In the case of $m=2$ this assertion is proven in
\cite{ZarhinTexel}.
\end{rem}

\begin{thm}
\label{mainL34} Let $K$ be a field of characteristic zero,
 $K_a$ its algebraic closure,
$f(x) \in K[x]$ an irreducible polynomial of degree $n$
 and
$\RR\subset K_a$ the set of its roots. We write $\Gal(f)$ for the
Galois group of $f$. Assume that there exist an odd power prime
$q\ge 3$ and an integer $m \ge 3$ such that $n=\frac{q^m-1}{q-1}$
and  the set $\RR$ could be identified with the projective space
$\P^{m-1}(\F_q)$ over the prime field $\F_q$ in such a way that
$\Gal(f)$ contains
 $\L_m(q):=\PSL_m(\F_{q})$ as a subgroup.
(E.g., $\Gal(f)=\L_m(q)$ or $\PGL_m(\F_{q})$.)
 Let $C_f$ be the
hyperelliptic curve $y^2=f(x)$. Let $J(C_f)$ be its jacobian,
$\End(J(C_f))$ the ring of $K_a$-endomorphisms of $J(C_f)$. Then
$\End(J(C_f))=\Z$.
\end{thm}

\begin{rem}
When $(m,q) \ne (4,3)$ this assertion is already proven is
\cite{ZarhinMMJ}. However, in this paper we give a unified proof
for all $(m,q)$.
\end{rem}

\section{Group theory}
\label{gt}
\begin{defn}
\label{ubd}
Suppose $\mathcal{G} \ne \{1\}$ is a perfect finite group and
$g>1$ is an integer. We say that $\mathcal{G}$ is $g$-{\sl
unbounded} if it enjoys the following properties:

\begin{enumerate}
\item[(i)]
Each homomorphism $\mathcal{G} \to \PSL(g-1,\C)$ is trivial;
\item[(ii)]
Either $g$ is odd and each homomorphism $\mathcal{G} \to
\PSL(g,\Q)=\SL(g,\Q)$ is trivial    or $g$ is even and each
homomorphism $\mathcal{G} \to \PSL(g,\R)=\SL(g,\R)/\{\pm 1\}$ is
trivial.
\end{enumerate}
\end{defn}

\begin{rem}
\label{unb}
\begin{enumerate}
\item[(i)]
Clearly, if every nontrivial irreducible projective representation
of  $\mathcal{G}$ in characteristic zero has dimension $>g$ then
all nontrivial projective representations of $\mathcal{G}$ in
characteristic zero have dimension $>g$ and therefore
$\mathcal{G}$ is $g$-unbounded.
\item[(ii)]
Clearly, if $\mathcal{H}$ is a simple non-abelian group isomorphic
to a subgroup of $\mathcal{G}$ then the $g$-unboundness of
$\mathcal{H}$ implies the $g$-unboundness of $\mathcal{G}$.
\item[(iii)]
Clearly, if every nontrivial  projective representation of
$\mathcal{G}$ in characteristic $2$ has dimension $>g$ then all
nontrivial projective representations of $\mathcal{G}$ in
characteristic zero have dimension $>g$ and therefore
$\mathcal{G}$ is $g$-unbounded.
\item[(iv)]
Every $\mathcal{G}$ is $2$-unbounded. Indeed, it suffices to check
that each homomorphism from $\mathcal{G}$ to $\PSL(2,\R)$ is
trivial. But each finite subgroup in $\PSL(2,\R)$ is the image of
a finite subgroup in $\SL(2,\R)$. Since each finite subgroup of
$\SL(2,\R)$ is commutative, the image $\pi(\mathcal{G})$ of each
homomorphism $\pi:\mathcal{G} \to \PSL(2,\R)$ is commutative. The
perfectness of $\mathcal{G}$ implies that $\pi(\mathcal{G})$ is
also perfect. This implies that $\pi(\mathcal{G})=\{1\}$, i.e.,
$\pi$ is trivial.
\end{enumerate}
\end{rem}


\begin{ex}
\label{bigA}
 Suppose $g\ge 4$ is an integer. Then the alternating groups $\A_{2g+1}$ and $\A_{2g+2}$ are $g$-bounded.
 Indeed, $2g+1>8$ and, by a theorem of Wagner \cite{Wagner}, every
 nontrivial linear representation of $\A_{2g+1}$ in characteristic
 $2$ has dimension $\ge 2g$. Since the Schur multiplier of
 $\A_{2g+1}$is $2$, one may easily deduce that every nontrivial
 projective representation of $\A_{2g+1}$ in characteristic
 $2$ also has dimension $\ge 2g>g$. By Remark \ref{unb}((iii),
 $\A_{2g+1}$is $g$-unbounded. Since
$\A_{2g+1}$ is isomorphic to a subgroup of $\A_{2g+2}$, the group
$\A_{2g+2}$ is also $g$-unbounded.
\end{ex}

\begin{exs}
\label{primer}
\begin{enumerate}
\item[(a)]
The Mathieu groups $\M_{23}$ and $M_{24}$ are $11$-unbounded.
Indeed, it follows from the Tables in \cite{Atlas} that all
nontrivial irreducible projective representations of $\M_{23}$ in
characteristic zero have dimension $\ge 22$. This implies that
$\M_{23}$ is $d$-unbounded for all $d<22$. Since $\M_{23}$ is
isomorphic to a subgroup of $\M_{24}$, the group $\M_{24}$ is also
$d$-unbounded for all $d<22$.
\item[(b)]
The Mathieu group $\M_{22}$  is $10$-unbounded. Indeed, it follows
from the Tables in \cite{Atlas} that in characteristic zero all
nontrivial irreducible projective representations of $\M_{22}$
have dimension $\ge 10$ and there are no nontrivial irreducible
{\sl linear} $10$-dimensional representations. This implies that
all homomorphisms $\M_{22}\to \PGL(9,\C)$ are trivial. Now in
order to establish the $10$-unboundness of $\M_{22}$ we need only
to check that each homomorphism from $\M_{22}$ to $\PSL(10,\R)$ is
trivial. Let as assume that $\rho: \M_{22} \to \PSL(10,\R)
=\SL(10,\R)/\{\pm 1\}$ is a nontrivial group homomorphism.
Clearly, $\rho$ lifts to a nontrivial homomorphism
$$\rho': \M_{22}' \to \SL(10,\R)$$
where $\M_{22}'$ sits in a central extension
$$\{1\} \to \{\pm 1\} \hookrightarrow \M_{22}' \twoheadrightarrow \M_{22} \to \{1\}.$$
This extension is non-splittable, since there are no nontrivial
homomorphisms from $\M_{22}$ to $\SL(10,\R)$ \cite{Atlas}.  On the
other hand, it also follows from the Tables in \cite{Atlas} that
there are no $10$-dimensional linear irreducible  representations
of $\M_{22}'$ defined over $\R$. This implies easily that there
are no nontrivial $10$-dimensional linear   representations of
$\M_{22}'$ defined over $\R$.This gives us the desired
contradiction and proves the $10$-unboundness of $\M_{22}$.
\item[(c)]
Suppose $r \ge 1$ and $m\ge 2$ are positive integers and assume
that $(m,r) \ne (2,1), (2,2), (3,2), (4,2)$. Let us put $q=2^r$.
We define a positive integer $g$ as follows.

$g=q/2$ if $m=2$ and $g=\frac{1}{2}\frac{q^{m}-q}{q-1}$ if $m>2$.

Then the group $\L_m(q)=\PSL_m(\F_{q})$ is $g$-unbounded. Indeed,
it is known \cite{TZ} that under our assumptions on $(m,q)$  all
nontrivial irreducible projective representations of $\L_m(q)$ in
characteristic zero have dimension $\ge \frac{q^{m}-q}{q-1}>g$ if
$m>2$. It is also known \cite{TZ} that all nontrivial irreducible
projective representations of $\L_2(q)$ in characteristic zero
have dimension $\ge q-1>g$ if $q>4$.
\item[(d)]
Suppose $q \ge 3$ is an odd power prime, $m\ge 3$ is an integer.
Let us put
$$g:=\left[\frac{1}{2} \left(\frac{q^m-1}{q-1} -1\right)\right]=
\left[\frac{1}{2} \left(\frac{q^m-q}{q-1}\right)\right].$$ Then
the group $\L_m(q)=\PSL_m(\F_q)$ is $g$-bounded. Indeed, it is
known \cite{TZ} that if $(q,m) \ne (3,4)$ then all nontrivial
irreducible projective representations of $\L_m(q)$ in
characteristic zero have dimension $\ge \frac{q^m-1}{q-1} -1>g$.

If $(q,m) = (3,4)$ then $g=19$. It follows from the Tables in
\cite{Atlas} that all nontrivial irreducible projective
representations of $\L_4(3)$ in characteristic zero have dimension
$\ge 26>19$.
\end{enumerate}
\end{exs}

\begin{ex}
\label{sm}
The Mathieu groups $\M_{11}$ and $M_{12}$ are $5$-unbounded.
Indeed, it follows from the Tables in \cite{Atlas} that all
nontrivial irreducible projective representations of $\M_{11}$ in
characteristic zero have dimension $\ge 10$. This implies that
$\M_{11}$ is $d$-unbounded for all $d<10$. Since $\M_{11}$ is
isomorphic to a subgroup of $\M_{12}$, the group $\M_{12}$ is also
$d$-unbounded for all $d<10$.
\end{ex}

\begin{defn}
\label{cvr}
Suppose $\mathcal{G}$ is a simple non-abelian group, $\ell$ is a
prime. Suppose
$$\{1\} \to N \hookrightarrow \mathcal{G}'
\stackrel{\gamma}{\twoheadrightarrow}
\mathcal{G}\to \{1\}$$ is a short exact sequence of finite groups
where $N$ is a group of exponent $1$ or $\ell$ and no proper
subgroup of $\mathcal{G}'$ maps {\sl onto} $\mathcal{G}$. (In
particular, $\mathcal{G}'$ is perfect.) Then $\mathcal{G}'$ is
called a {\sl minimal} {\sl cover} of $\mathcal{G}$.

In addition, if either $N=\{1\}$ (i.e.,
$\mathcal{G}'=\mathcal{G}$) or the exponent of $N$ is $\ell$ then
we say that $\mathcal{G}'$ is a {\sl minimal} $\ell$-{\sl cover}
of $\mathcal{G}$.
\end{defn}

\begin{rem}
\label{cover}
 Clearly, the minimal cover $\mathcal{G}'$ is always
perfect. It is also clear that each normal subgroup in
$\mathcal{G}'$ except $\mathcal{G}'$ itself lies in $N$. This
implies easily that if $\rho: \mathcal{G}' \to M$ is a nontrivial
group homomorphism then $\ker(\rho)$ lies in $N$ and therefore the
image $\rho(\mathcal{G}')$ is also a minimal cover of
$\mathcal{G}$. In addition, if $\mathcal{G}'$ is a minimal
$\ell$-cover then $\rho(\mathcal{G}')$ is also one.
\end{rem}

\begin{rem}
\label{sur} Suppose $\mathcal{G}$ is a simple non-abelian group
and $\gamma:H \twoheadrightarrow \mathcal{G}$ is a surjective
homomorphism of finite groups. Let $\mathcal{G}'$ be a subroup of smallest order among the
 subgroups $H'$ of $H$ such that
$\gamma(H')=\mathcal{G}$. (Clearly, such $\mathcal{G}'$ always
exists.) Then $\gamma: \mathcal{G}'\to \mathcal{G}$ is a minimal
cover. In addition, if the kernel of is either trivial or has
exponent $\ell$ then $\gamma: \mathcal{G}'\to \mathcal{G}$ is a
minimal $\ell$-cover.
\end{rem}

\begin{exs}
\label{seven}
\begin{enumerate}
\item[(i)]
Suppose $\mathcal{G}$ is a simple non-abelian group isomorphic
either to $\A_5\cong \L_2(4)$ or to $\A_6$. Suppose $\mathcal{G}'$
is a minimal cover of  $\mathcal{G}$.
 Then, by Remark \ref{unb}(iv), the perfect group $\mathcal{G}'$ is  $2$-bounded.
\item[(ii)]
Suppose $\mathcal{G}$ is a simple non-abelian group isomorphic
either to $\L_3(2)\cong \L_2(7)$ or to $\A_7$ or to
$\A_8\cong\L_4(2)$. Notice that in all these cases the order of
$\mathcal{G}$ is divisible by $7$.

 Suppose $\mathcal{G}'$ is a minimal cover of  $\mathcal{G}$.
 Then $\mathcal{G}'$ is  $3$-bounded. Indeed, if
$\pi_2:\mathcal{G}' \to \PSL(2,\C)$ is a {\sl nontrivial}
homomorphism then, by Remark \ref{cover},  its image
$H':=\pi(\mathcal{G}')\subset \PSL(2,\C)$ is also a minimal cover
of  $\mathcal{G}$. In particular, $H_2$ is a perfect finite group
having a quotient isomorphic to  $\mathcal{G}$. This implies that
nonsolvable $H_2\subset \PSL(2,\C)$ is not isomorphic to $\A_5$
which could not be true (\cite{Suzuki}, Th. 6.17 on p. 404). The
obtained contradiction implies that there are no nontrivial
homomorphisms from $\mathcal{G}'$ to $\PSL(2,\C)$.

Now assume that there exists a {\sl nontrivial} homomorphism
$\pi_3:\mathcal{G}' \to \SL(3,\Q)$. As above, the image
$H_3:=\pi(\mathcal{G}')\subset \SL(3,\Q)$ is also a minimal cover
of  $\mathcal{G}$. In particular, $H_3$ is a  finite group having
a quotient isomorphic to  $\mathcal{G}$; in particular, $7$
divides the order of $H_3$ and therefore $H_3 \subset \SL(3,\Q)$
contains an element of order $7$. But this is not true, since the
degree of the $7$th cyclotomic field over $\Q$ is $6>3$. The
obtained contradiction ends the proof of the $3$-unboundness of
$\mathcal{G}'$.
\end{enumerate}
\end{exs}

 We will use the following result of Feit-Tits
(\cite{FT}, Theorem on pp. 1092--1093 and Prop. 4.1 on p. 1098)
concerning complex projective representations of minimum degree.
(See also \cite{KL}.)

\begin{thm}[Feit-Tits Theorem]
\label{FTKL}
 Suppose $\mathcal{G}$ is a known simple non-abelian group that
 is not a a group of Lie type in characteristic $2$.
 Suppose $\mathcal{G}' \stackrel{\gamma}{\twoheadrightarrow}
\mathcal{G}$ is a minimal cover of $\mathcal{G}$ and $d$ is the
smallest positive integer such that there exists a nontrivial
homomorphism $\mathcal{G}' \to \PGL(d,\C)$.
 Then the kernel of each homomorphism $\mathcal{G}' \to \PGL(d,\C)$ contains
$\ker({\gamma})$. In particular,  $\mathcal{G}$ is isomorphic to a
subgroup of $\PGL(d,\C)$.
\end{thm}

\begin{rem}
\label{Lmq0}
 Kleidman and Liebeck \cite{KL} studied the case of
simple groups of Lie type in characteristic $2$. In particular,
they proved the following assertion (ibid, Th. 3 on p. 316).
Suppose  $q$ is a power of $2$ and $m \ge 2$ is an integer such
that $(m,q)) \ne (2,2)$ (i.e., $\L_m(q)$ is a simple non-abelian
group). If $g$ is a positive integer and $\mathcal{G}'$ is a
minimal cover of $\L_m(q)$ such that $\mathcal{G}'$ is isomorphic
to a subgroup of $\PGL(g,\C)$ then either $m=4$ and $g \ge q^3$ or
$\L_m(q)$ is isomorphic to  a subgroup of $\PGL(g,\C)$.
\end{rem}

\begin{lem}
\label{Lmq} Suppose $q \ge 2$ is an integral power of $2$ and $m
\ge 2$ is an integer such that $(m,q) \ne (2,2), (2,4), (3,2),
(4,2)$. Suppose $\mathcal{G}'$ is a minimal cover of $\L_m(q)$.
\begin{enumerate}
\item[(i)]
If $m>2$ then $\mathcal{G}'$ is
$\frac{1}{2}\frac{q^m-q}{q-1}$-unbounded.
\item[(ii)]
If $m>2$ then $\mathcal{G}'$ is $\frac{q}{2}$ -unbounded.
\end{enumerate}
\end{lem}

\begin{proof}
Let us start with the case of $m=2$. Let us put $g=\frac{q}{2}$.
Let $\rho:\mathcal{G}'\to \PGL(g,\C)$ is a {\sl nontrivial} group
homomorphism. By Remark \ref{cover}, the image
$\rho(\mathcal{G}')\subset\PGL(g,\C)$ is also a minimal cover of
$\L_2(q)$. Applying Remark \ref{Lmq0} to $\rho(\mathcal{G}')$, we
conclude that $\L_2(q)$ is isomorphic to a subgroup of
$\PGL(g,\C)$. By Example \ref{primer}(c), this is not true. The
obtained contradiction implies that there are no nontrivial
homomorphisms from $\mathcal{G}'$ to $\PGL(g,\C)$.

Now assume that $m>2$. Let us put
$g=\frac{1}{2}\frac{q^m-q}{q-1}$. Notice that if $m=4$ then
$$g=\frac{1}{2}(q^3+q^2+q)<q^3.$$
 Let $\rho:\mathcal{G}'\to \PGL(g,\C)$ is a {\sl
nontrivial} group homomorphism. Again, the image
$\rho(\mathcal{G}')\subset\PGL(g,\C)$ is also a minimal cover of
$\L_m(q)$. Applying Remark \ref{Lmq0} to $\rho(\mathcal{G}')$, we
conclude that either $m=4$ and $g \ge q^3$ or $\L_m(q)$ is
isomorphic to a subgroup of $\PGL(g,\C)$. But we have already seen
that if $m=4$ then $g<q^3$. This implies that $\L_m(q)$ is
isomorphic to a subgroup of $\PGL(g,\C)$. By Example
\ref{primer}(c), this is not true. The obtained contradiction
implies that there are no nontrivial homomorphisms from
$\mathcal{G}'$ to $\PGL(g,\C)$.
\end{proof}

\begin{lem}
\label{primer1} Suppose a simple non-abelian finite group
$\mathcal{G}$ and an integer $g>1$ enjoy one of the following
properties:
\begin{enumerate}
\item[(a)]
$g=11$ and $\mathcal{G}=\M_{23}$ or $\M_{24}$;
\item[(b)]
$g=10$ and $\mathcal{G}=\M_{22}$;
\item[(c)]
$g=5$ and $\mathcal{G}=\M_{11}$ or $\M_{12}$;
\item[(d)]
$g=\left[\frac{1}{2} \left(\frac{q^m-1}{q-1} -1\right)\right]$
where $q \ge 3$ is an odd power prime, $m \ge 3$ is an integer and
$\mathcal{G}=\L_m(q)$;
\item[(e)]
$\mathcal{G}=\A_{2g+1}$ or $\A_{2g+2}$;
\end{enumerate}
If $\mathcal{G}'$ is a minimal cover of $\mathcal{G}$ then
$\mathcal{G}'$ is $g$-unbounded.
\end{lem}

\begin{proof}

Case (a). We have seen that there are no nontrivial homomorphisms
to $\PSL(11,\C)$ either from $\M_{23}$ or $\M_{24}$. It follows
from the Feit--Tits theorem that the same is true for the minimal
cover $\mathcal{G}'$.

Case (c). We have seen that there are no nontrivial homomorphisms
to $\PSL(5,\C)$ either from $\M_{11}$ or $\M_{12}$. It follows
from the Feit--Tits theorem that the same is true for the minimal
cover $\mathcal{G}'$.

Case (d).  It follows from the Feit--Tits theorem combined with
Example \ref{primer}(d) that there are no nontrivial homomorphisms
from $\mathcal{G}'$ to $\PSL(g,\C)$.


Case (b). Since every homomorphism from $\M_{22}$ to $\PSL(9,\C)$
is trivial, we conclude that, thanks to the Feit--Tits theorem,
that every homomorphism from $\mathcal{G}'$  to $\PSL(9,\C)$ is
also trivial. In order to finish the proof we have to check that
every homomorphism from $\mathcal{G}'$  to $\PSL(10,\R)$ is also
trivial. Let us assume that $\rho': \mathcal{G}' \to \PSL(10,\R)
=\SL(10,\R)/\{\pm 1\}$ is a nontrivial group homomorphism. Let us
consider the composition
$$\pi:\mathcal{G}' \stackrel{\rho'}{\longrightarrow} \PSL(10,\R) \subset \PSL(10,\C)$$
of $\rho'$ and the natural embedding $\PSL(10,\R) \subset
\PSL(10,\C)$. Clearly, the composition
$\pi:\mathcal{G}'\to\PSL(10,\C)$ is a {\sl nontrivial} group
homomorphism. Since $\mathcal{G}'$ is a minimal cover of $\M_{22}$
and $10$ is the smallest dimension of a nontrivial projective
representation of $\M_{22}$ over $\C$ \cite{Atlas}, the Feit--Tits
theorem implies that $\ker(\pi)$ contains
$\ker(\mathcal{G}'\twoheadrightarrow \M_{22})$. Since the image of
$\pi$ lies in $\PSL(10,\R) \subset \PSL(10,\C)$, we conclude that
$\pi$ gives rise to a nontrivial homomorphism $\M_{22} \to
\PSL(10,\R)$. Contradiction.

Case (e) follows easily from the Feit--Tits theorem combined with
Examples  \ref{seven} and Lemma \ref{Lmq}.
\end{proof}

\begin{thm}
\label{main0} Suppose $g>1$ is an integer, $D$ a
finite-dimensional semisimple $\Q$-algebra enjoying the following
properties:
\begin{enumerate}
\item[(i)]
Let us split $D$ into a direct sum
$$D =\oplus_{i=1}^r D_i$$
of simple $\Q$-algebras $D_i$. Then the number $r$ of summands
does not exceed $g$;
\item[(ii)]
Let us present a summand $D_i$ as the algebra of square matrices
of size $d_i$ over a division $\Q$-algebra $T_i$. Then every $T_i$
admits a positive involution. In addition,
$$\oplus_{i=1}^r d_i \le g.$$

\item[(iii)]
If $r=1$ (i.e., $D=D_1$ is simple) then $n_1 \cdot \dim_{\Q}(T_1)$
divides $2g$. In addition, the center of $D_1$ is either a totally
real number field of degree dividing $g$ or a CM-field of degree
dividing $2g$. Also, if $n_1 \cdot \dim_{\Q}(T_1)=2g$ and $T_1$ is
a quaternion $\Q$-algebra then it is indefinite.
\end{enumerate}

 Suppose $H$ is a $g$-unbounded group and $\rho: H \to
\Aut(D)$ is a group homomorphism such that the subalgebra
$$D^H=\{u \in D\mid \rho(h)u=u \quad \forall h \in H\}$$ of
$H$-invariants coincides with $\Q$. Then $D=\Q$.
\end{thm}

\begin{proof}
Let $C_i$ be the center of $D_i$. Then $C_i$ is either a totally
real number field or a CM-field. Clearly,
$$C=\oplus _{i=1}^r C_i$$
is the center of $D$. This implies that $C$ is $H$-stable and the
action of $H$ permutes $C_i$'s. This gives rise to a homomorphism
from $H$ to the group $\SS_r$ of permutations in $r$ letters which
must be trivial. Indeed, the perfectness of $H$ implies that its
image in $\SS_r$ lies in the alternating subgroup $\A_r$  which
embeds into $\PSL(r-1,\C)$ if $r>2$ and the inequality $r \le g$
and the triviality of homomorphisms implies in this case that $H
\to \SS_r$ is trivial. If $r \le 2$ then $\A_r$ is itself trivial.
So, $H$ leaves stable each $C_i$. This implies easily that
$\oplus_{i=1}^r\Q$ consists of $H$-invariants. Since $D^H=\Q$, we
conclude that $r=1$ and therefore $D=D_1$ and $C=C_1$ is also the
center of $M_1$. So, $C$ is either a totally real number field of
degree dividing $g$ or a purely imaginary quadratic extension of a
totally real number field $C^{+}$ where $[C^{+}:\Q]$ divides $g$.
In the case of totally real $C$ let us put $C^{+}:=C$. Clearly, in
both cases $C^{+}$ is the largest totally real subfield of $C$ and
therefore the action of $H$ leaves $C^{+}$ stable. Let us put
$d:=[C^{+}:\Q]$. I claim that $d=1$, i.e., $C^{+}=\Q$. Indeed,
suppose $d>1$. Clearly, one may identify $\Aut(H)$ with a subgroup
of $\GL(d-1,\Q)$ and therefore the action of $H$ on $C^{+}$ gives
us a homomorphism
$$H \to \Aut(C^{+}) \subset \GL(d-1,\Q),$$
whose triviality we need to check. Assume the contrary. The
perfectness of $H$ implies that the {\sl nontrivial} image of the
composition $H \to \GL(d-1,\Q)$ is, in fact, a perfect subgroup of
$\SL(d-1,\Q)$. This perfectness implies, in turn, that the image
of $H$ in $\PSL(d-1,\Q)$ is also nontrivial. Taking into account
the inequality  $d \le g$ and the inclusion $\PSL(d-1,\Q)\subset
\PSL(d-1,\C)$, we obtain a nontrivial homomorphism $H \to
\PSL(d-1,\C)$. This contradicts to the $g$-unboundness of $H$.
Hence $d=1$ and $C^{+}=\Q$.

 Now I claim that $C=\Q$. Indeed, if $C \ne \Q$ then $C$ is an
imaginary quadratic field and $\Aut(C)$ is a cyclic group of order
$2$. The perfectness of $H$ implies that $H \to \Aut(C)$ is
trivial and therefore $C$ consists of $H$-invariants. Since
$D^H=\Q$, we get a contradiction. Hence $C=\Q$.

So, $T_1$ is a central simple $\Q$-algebra with a positive
involution. Hence either $T_1=\Q$ or  a quaternion $\Q$-algebra.

Assume that $\T_1=\Q$. Then $D=D_1$ is the matrix algebra of size
$d_1$ over $\Q$. Clearly, $d_1 \le g$. By Skolem-Noether theorem,
$\Aut(D)$ is $\PGL(d_1,\Q)$. So, perfect $H$ acts on $D$ via
$\rho:H \to\Aut(D)=\PGL(d_1,\Q)$, whose image must lie in
$\PSL(d_1,\Q)$. Since $d_1 \le g$, $\PSL(d_1,\Q)$ is a subgroup of
$\PSL(g,\Q)$, we obtain the triviality of $\rho:H \to \Aut(D)=
\PGL(d_1,\Q)\subset \PSL(g,\Q)$. This implies that the whole $D$
consists of $H$-invariants. Since $D^H=\Q$, we conclude that
$D=\Q$.

Now assume that $T_1$ is a quaternion $\Q$-algebra. Then $D \ne
\Q$ and therefore $\rho: H \to \Aut(D)$ is nontrivial. We need to
arrive to a contradiction.

We have  $\dim_{\Q}(D_1)=4$ and $4 n_1 \le 2g$.

Assume that $4 n_1=2g$. Then $g= 2 n_1$ is {\sl even} and
$T_1\otimes_{\Q}\R$ is isomorphic to the matrix algebra of size
$2$ over $\R$. This implies that $D_{\R}:=D\otimes_{\Q}\R$ is the
matrix algebra of size $g$ and
$$\Aut(D) \subset \Aut_{\R}(D_{\R})=\PGL(g,\R).$$
Therefore the nontrivial $\rho$ gives rise to a nontrivial
homomorphism $H \to \PGL(g,\R)$. Again, the perfectness of $H$
implies that the image lies in $\PSL(g,\R)$ and we get a
nontrivial homomorphism $H \to \PSL(g,\R)$. Contradiction.

Assume that $4 n_1<2g$. Then $D_{\C}:=D\otimes_{\Q}\C$ is the
matrix algebra over $\C$ of size $2n_1<g$ and
$$\Aut(D) \subset \Aut_{\C}(D_{\C})=\PGL(2n_1,\C).$$
Therefore the nontrivial $\rho$ gives rise to a nontrivial
homomorphism $H \to \PGL(g-1,\C)$. Again, the perfectness of $H$
implies that the image of $H$ lies in $\PSL(g-1,\C)$ and we get a
nontrivial homomorphism $H \to \PSL(g-1,\C)$. Contradiction.
\end{proof}

Let $B$ be a finite set consisting of $n \ge 5$ elements. We write $\Perm(B)$ for the group of permutations of $B$. A choice of ordering on $B$ gives rise to an isomorphism
$$\Perm(B) \cong \Sn.$$
 Let us consider the permutation module $\F_2^{B}$: the
$\F_2$-vector space of all functions $\varphi:B \to \F_2$. The
space $\F_2^{B}$ carries a natural structure of
$\Perm(B)$-module and contains the stable line $\F_2\cdot 1_{B}$
of constant functions and the stable hyperplane $(\F_2^{B})^0$
of functions $\varphi$ with $\sum_{\alpha\in
B}\varphi(\alpha)=0$. Clearly, $(\F_2^{B})^0$ contains
$\F_2\cdot 1_{B}$ if and only if $n$ is even. Let us put
$Q_{B}:=(\F_2^{B})^0$ if $n$ is odd and
$Q_{B}:=(\F_2^{B})^0/(\F_2\cdot 1_{B})$ if $n$ is even.
Clearly, $Q_{B}$ carries a natural structure of faithful
$\Perm(B)$-module. For each permutation group $H\subset
\Perm(B)$ the corresponding $H$-module is called the {\sl heart} of
the permutation representation of $H$ on $B$ over $\F_2$
\cite{Klemm}, \cite{Mortimer}, \cite{Ivanov}.

\begin{lem}
\label{Kl}
$\End_H(Q_{B})=\F_2$ if either $n$ is odd and
$H$ acts $2$-transitively on $B$ or $n$ is even and $H$ acts
$3$-transitively on $B$.
\end{lem}

\begin{proof}
See Satz 4 in \cite{Klemm}.
\end{proof}

\begin{lem}
\label{qodd}
Suppose $q>2$ is an integral {\sl odd} power prime, $m \ge 3$ is an integer,
$B=\P^{m-1}(\F_q)$ is the corresponding $(m-1)$-dimensional projective space over the finite field $\F_q$
 and $H=\L_m(q)=\PSL_m(\F_q) \subset \Perm(B)$
is the corresponding projective special linear group over $\F_q$
acting naturally and
 faithfully on the projective space. Then the $H$-module $Q_B$ ia absolutely simple.
 In particular, $\End_H(Q_{B})=\F_2$.
\end{lem}

\begin{proof}
See \cite{Mortimer}, Table 1 on page 2.
\end{proof}

\begin{rem}
\label{Lmq2}
 Suppose $q = 2^r$ is an integral power of $2$ and $m\ge 2$ is an integer
such that either $m \ge 3$ or $q>4$. Suppose $B=\P^{m-1}(\F_q)$ is
 the corresponding
$(m-1)$-dimensional projective space over the finite field $\F_q$
 and $H=\L_m(q)=\PSL_m(\F_q) \subset \Perm(B)$
is the corresponding projective special linear group over $\F_q$
acting naturally and faithfully on the projective space. Then
$\End_H(Q_{B})=\F_2$. Indeed,  it is well-known that
$H=\PSL_m(\F_q)$ acts doubly transitively on $B=\P^{m-1}(\F_q)$.
Clearly, $\#(B)=\frac{q^m-1}{q-1}$ is odd and therefore the
assertion follows from Lemma \ref{Kl}. Notice that if $m \ge 3$
then the $H$-module $Q_B$ is reducible \cite{Mortimer}; see \S 5
of \cite{Ivanov} for details. If $m=2$ then the $H$-module $Q_B$
is absolutely simple \cite{Mortimer} (and even very simple
\cite{ZarhinTexel}).
\end{rem}

\section{Endomorphisms of abelian varieties}
\label{AV} Let  $K$ be a field of characteristic zero. We fix its algebraic closure $K_a$
and write $\Gal(K)$ for the absolute Galois group $\Aut(K_a/K)$.
Let $X$ is an abelian variety of positive dimension defined over $K$. Then the group
$X(K_a)$ of its algebraic points has a natural structure of
$\Gal(K)$-module. If $m$ is a positive integer  then we write
$X_m$ for the kernel of multiplication by $m$ in $X(K_a)$. It is
well-known \cite{MumfordAV} that $X_m$ is a free $\Z/m\Z$-module
of rank $2\dim(X)$ provided with the structure of $\Gal(K)$-module
inherited from $X(K_a)$. We denote by $\tilde{\rho}_{m,X}$  the
corresponding homomomorphism
$$\tilde{\rho}_{m,X}: \Gal(K) \to \Aut(X_{m}) \cong
\GL(2\dim(X),\Z/m\Z)$$ which defines the structure of Galois
module on $X_{m}$. We have
$$\tilde{\rho}_{m,X}(\sigma)(x)=\sigma(x) \quad \forall \sigma
\in \Gal(K), x \in X_{m} \subset X(K_a).$$ We write
$$\tilde{G}_{m,X,K}=\tilde{\rho}_{m,X}(\Gal(K))$$
for  the image of $\Gal(K)$ in $\Aut(X_{m})$. If $K(X_m)$ is the
field of definition of all points of order $m$ on $X$ then it is a
finite Galois extension of $K$, whose Galois group
$\Gal(K(X_m)/K)=\tilde{G}_{m,X,K}$.

Suppose $m>1$. Clearly, $X_m= m X_{m^2}$ coincides with the kernel
of multiplication by $m$ in $X_{m^2}$. In particular, every
endomorphism of the commutative group $X_{m^2}$ leaves $X_m$
stable. Therefore the restriction to $X_m$ gives rise to a natural
(obviously surjective) ring homomorphism (the reduction modulo
$m$)

$$\red_m:\End(X_{m^2}) \twoheadrightarrow  \End(X_{m}); \quad \red_m(u)(x)=ux \quad \forall u \in \End(X_{m^2}),
x \in  X_{m}.$$ Clearly, $\ker(\red_m)=m\End(X_{m^2})$. This
implies that each $v \in \ker(\red_m)$ satisfies $v^2=0=mv$.

Restricting $\red_m$ to the automorphism group $\Aut((X_{m^2})$ of
$X_{m^2}$, we obtain
 the (obviously surjective) group homomorphism
$$\red_m^*:\Aut(X_{m^2}) \twoheadrightarrow  \Aut(X_{m});\quad \red_m^*(u)(x)=ux \quad \forall u \in \Aut(X_{m^2}),
x \in  X_{m}.$$ Clearly, $\ker(\red_m^*)=\I+m\End(X_{m^2})$. (Here
$\I$ is the identity automorphism of $X_{m^2}$.) This implies that
each $u \in \ker(\red_m^*)$ is of the form $\I+v$ with $v^2=0=mv$.
This implies that $u^m=\I$.

\begin{rem}
Notice that the homomorphisms $\red_m$ and $\red_m^*$ do not
depend on the choice of the field of definition $K$ for $X$. In
particular, they both are $\Gal(K)$-equivariant.
\end{rem}

 Clearly,
$$\tilde{\rho}_{m^2,X}=\red_m^*\tilde{\rho}_{m,X};$$
in particular,
$$\tilde{G}_{m,X,K}=\red_m^*(\tilde{G}_{m^2,X,K}).$$

 There is an important special case when $m=\ell$ is a
prime. Then $X_m=X_{\ell}$ is a $2\dim(X)$-dimensional
$\F_{\ell}$-vector space provides with the structure of
$\Gal(K)$-module inherited from $X(K_a)$. We have $u^{\ell}=\I$
for each $u \in \ker(\red_{\ell}^*)$ and therefore the exponent of
the nontrivial finite group $\ker(\red_{\ell}^*)$ is $\ell$. We
also have
$$\tilde{G}_{\ell,X,K}=\red_{\ell}^*(\tilde{G}_{\ell^2,X,K})$$
and therefore the kernel of the surjective group homomorphism
$$\red_{\ell}^*:\tilde{G}_{\ell^2,X,K} \twoheadrightarrow
\tilde{G}_{\ell,X,K}$$ is either trivial or a finite group of
exponent $\ell$.

\begin{rem}
\label{ell2} Assume that $\tilde{G}_{\ell,X,K}$ contains a simple
non-abelian subgroup $\mathcal{G}$. Let $H \subset
\tilde{G}_{\ell^2,X,K}$ be the preimage of $\mathcal{G}$ with
respect to $\red_{\ell}^*:\tilde{G}_{\ell^2,X,K}
\twoheadrightarrow \tilde{G}_{\ell,X,K}$. Clearly,
$\red_{\ell}^*:H \to \mathcal{G}$ is a surjective homomorphism,
whose kernel is either trivial or has exponent $\ell$. According
to Remark \ref{sur} there exists a subgroup
$$\mathcal{G}' \subset H \subset \tilde{G}_{\ell^2,X,K}$$
such that $\red_{\ell}^*:\mathcal{G}' \twoheadrightarrow
\mathcal{G}$ is a minimal $\ell$-cover.

Clearly, $\{\sigma\in\Gal(K)\mid \rho_{\ell^2,X}(\sigma) \in
\mathcal{G}'\}$ is an open  subgroup of finite index in $\Gal(K)$
and therefore coincides with the Galois group  of certain finite
separable algebraic extension $L=L_{\mathcal{G}}$ of $K$. Clearly,
we have
$$\rho_{\ell^2,X}(\Gal(L))=\mathcal{G}', \quad
\rho_{\ell,X}(\Gal(L))=\red_{\ell}^*(\mathcal{G}')=\mathcal{G}.$$
In other words,
$$\tilde{G}_{\ell^2,X,L}=\mathcal{G}',\quad
\tilde{G}_{\ell,X,L}=\mathcal{G}$$ and
$\red_{\ell}^*:\tilde{G}_{\ell^2,X,L}\twoheadrightarrow
\tilde{G}_{\ell,X,L}$ is a minimal $\ell$-cover.
\end{rem}

We write $\End(X)$ for the ring of all $K_a$-endomorphisms of $X$
and $\End_K(X)$ for the ring of all $K$-endomorphisms of $X$. We
have

$$\Z=\Z\cdot \I_X \subset \End_K(X) \subset \End(X)$$ where $\I_X$
is the identity automorphism of $X$.

Since $X$ is defined over $K$, one may associate with every $u \in
\End(X))$ and $ \sigma \in \Gal(K)$ an endomorphism $^{\sigma}u
\in \End(X)$ such that

$$^{\sigma}u(x)=\sigma u(\sigma^{-1}x)
\quad \forall x \in X(K_a).$$ In fact, there is a group
homomorphism

$$\kappa_{X}: \Gal(K) \to \Aut(\End(X)); \quad
\kappa_{X}(\sigma)(u)=^{\sigma}u \quad \forall \sigma \in \Gal(K),
u \in \End(X).$$

It is well-known that $\End_K(X)$ coincides with the subring of
$\Gal(K)$-invariants in $\End(X)$, i.e., $$\End_K(X)=\{u\in
\End(X)\mid ^{\sigma}u=u \quad \forall \sigma \in \Gal(K)\}.$$ It
is also well-known, that $\End(X)$, viewed as a group (with
respect to addition) is a free commutative group of finite rank
and $\End_K(X)$ is its {\sl pure} subgroup, i.e., the
quotient-group $\End(X)/\End_K(X)$ is also  a free commutative
group of finite rank. It is also well-known that there exists a
finite Galois extension $K'$ of $K$ such that all the
endomorphisms of $X$ are defined over $K'$, i.e.
$$\Gal(K') \subset \ker(\kappa_{X})\subset \Gal(K).$$

\begin{rem}
\label{silver}
 It is proven in \cite{Silverberg} that all the endomorphisms of
$X$ are defined over $K(X_4)\bigcap K(X_{\ell})$ where $\ell$ is
an arbitrary odd prime. This means that
$\ker(\tilde{\rho}_{4,X})\subset \ker(\kappa_{X})$ and
$\ker(\tilde{\rho}_{\ell,X})\subset \ker(\kappa_{X})$ for all {\sl
odd} primes $\ell$. This implies that if $\Gamma_K
=\kappa_{X}(\Gal(K)) \subset \Aut(\End(X))$ then there exists a
surjective homomorphism
$$\tilde{G}_{4,X} \twoheadrightarrow \Gamma_K\subset \Aut(\End(X))$$ and surjective
homomorphisms
$$\tilde{G}_{\ell,X} \twoheadrightarrow \Gamma_K\subset \Aut(\End(X))$$
for all odd primes $\ell$. Notice that
$$\End_K(X)=\End(X)^{\Gamma_K}.$$
\end{rem}

\begin{rem}
\label{endo}
\begin{enumerate}
\item[(i)]
Let us put $\End^0(X):=\End(X)\otimes\Q$. It is well-known
(\cite{MumfordAV}, \S 21) that $\End^0(X)$ is a semisimple
finite-dimensional $\Q$-algebra. Clearly, the natural map
$$\Aut(\End(X)) \to \Aut(\End^0(X))$$
is an embedding.
\item[(ii)]
Recall that $X$ is isogenous over $K_a$ to a product
$\prod_{i=1}^r Y_i^{d_i}$ where $Y_i$'s are mutually non-isogenous
absolutely simple abelian varieties (of positive dimension) over
$K_a$. Then
$$g=\dim(X)=\sum_{i=1}^r d_i\cdot\dim(Y_i) \ge \sum_{i=1}^r d_i.$$
Let us put $T_i:=\End^0(Y_i)$ and denote by $D_i$ the  algebra of
square matrices of size $d_i$ over $T_i$. Then each $T_i$ is a
division $\Q$-algebra admitting a positive involution
(\cite{MumfordAV}, \S 21).   Let us denote by $C_i$  the center of
$T_i$. Then either $C_i$ is a totally real number field and
$[C_i:\Q]$ divides $\dim(Y_i)$ or $C_i$ is a CM-field and
$[C_i:\Q]$ divides $2\dim(Y_i)$. It is also clear that $C_i$ is
the center of $D_i$.

 Since $\fchar(K_a)=0$, the number
$\dim_{\Q}(T_i)$ divides $2\dim(Y_i)$ (\cite{MumfordAV}, \S 21, p. 202). Clearly,
$D_i=\End^0(Y_i^{d_i})$ and
$$\End^0(X)=\oplus_{i=1}^r D_i.$$

\item[(iii)]
Assume now that $r=1$, i.e., $X$ is isogenous to $Y_1^{d_1}$ and
$\End^0(X)=D_1$. Then $g=\dim(X)=d_1\cdot \dim(Y_1)$. Hence either $C_1$
is totally real number field and $[C_1:\Q]$ divides $g$ or $C_i$
is a CM-field and $[C_i:\Q]$ divides $2g$. It is also clear that
$d_1\dim_{\Q}(T_1)$ divides $d_1\cdot 2\dim(Y_1)=2\dim(X)=2g$. If
$C_1=\Q$ then $T_1$ is either $\Q$ or a quaternion $\Q$-algebra.

\item[(iv)]
We continue to assume that $r=1$. If $T_1$ a quaternion
$\Q$-algebra and $d_1\dim_{\Q}(T_1)=2g$ then, taking into account
that $\dim_{\Q}(T_1)=4$, we conclude that $g$ is even, $d_1=g/2$
and $Y_1$ is an absolutely simple abelian surface. Since in
characteristic zero the endomorphism algebra of an absolutely
simple abelian surface  is either a field or an indefinite
quaternion $\Q$-algebra \cite{Oort} (see also \cite{OZ}), we
conclude that $T_1=\End^0(Y_1)$ is an {\sl indefinite} quaternion
$\Q$-algebra.
\end{enumerate}
\end{rem}

\begin{thm}
\label{apps} Suppose   $K$ is a field of characteristic $0$,
suppose $X$ is an abelian variety over a $K$ of dimension $g>1$.
 Suppose $\ell$ is a prime,
$$\tilde{G}_{\ell,X,K}=\tilde{\rho}_{\ell,X}(\Gal(K))\subset \Aut(X_{\ell})$$
is the image of $\Gal(K)$ in $\Aut(X_{\ell})$. Let us put
$g=\dim(X)$ and assume that $g>1$ (i.e., $X$ is not an elliptic
curve). Assume that $\tilde{G}_{\ell,X,K}$ contains a simple
non-abelian subgroup $\mathcal{G}$ such that
$$\End_{\mathcal{G}}(X_{\ell})=\F_{\ell}$$ and one of the
following conditions holds:
\begin{enumerate}
\item[(a)]
$\ell$ is odd and $\mathcal{G}$ is $g$-unbounded.
\item[(b)]
$\ell=2$ and every $2$-minimal cover of $\mathcal{G}$ is
$g$-unbounded.
\end{enumerate}
 Then the ring $\End(X)$ of all $K_a$-endomorphisms of $X$
coincides with $\Z$.
\end{thm}

\begin{proof}
First, using Remark \ref{ell2} we may replace $K$ by its finite
separable algebraic extension $L$ in such a way that
$\tilde{G}_{\ell,X,L}=\mathcal{G}$ and
$$\red_{\ell}^*:\tilde{G}_{\ell^2,X,L}\hookrightarrow
\tilde{G}_{\ell,X,L}=\mathcal{G}$$ is a minimal $\ell$-cover.
Clearly, $\tilde{G}_{\ell,X,L}$ is $g$-unbounded if $\ell$ is odd.
If $\ell=2$ then it follows from Remark \ref{ell2} that $\tilde{G}_{\ell,X,L}$
is a minimal $2$-cover of $\tilde{G}_{2,X,L}=\mathcal{G}$ and therefore is $g$-unbounded.

 Second, I claim that $\End_L(X)=\Z$. Indeed, it is
well-known that there is an embedding.
$$\End_L(X)\otimes \Z/\ell\Z \hookrightarrow
\End_{\Gal(K)}(X_{\ell}).$$
On the other hand, since
$\rho_{\ell,X}(\Gal(L))=\tilde{G}_{\ell,X,L}$,
$$\End_{\Gal(K)}(X_{\ell})=\End_{\tilde{G}_{\ell,X,L}}(X_{\ell})=
\End_{\mathcal{G}}(X_{\ell})=\F_{\ell},$$ the rank of free
commutative group $\End_L(X)$ is either $0$ or $1$. Clearly, it
must be $1$ and this implies that $\End_L{X}=\Z$.

Now let us put $D:=\End^0(X)$. Clearly, $\End(X)$ is a
$\Z$-lattice in the $\Q$-vector space $D$. Let $\Aut(D)$ be the
group of automorphisms of the $\Q$-algebra $D$. We have
$\Aut(\End(X)) \subset \Aut(D)$. We have
$$\kappa_X(\Gal(L))=\Gamma_L \subset\Aut(\End(X)) \subset
\Aut(D).$$ Clearly, we have
$$D^{\Gamma_L}=\End(X)^{\Gamma_L}\otimes\Q=\End_L(X)\otimes\Q=\Z\otimes\Q=\Q.$$

We are going to finish the proof, using Theorem \ref{main0}. Let us put
$H:=\tilde{G}_{\ell,X,L}$ if $\ell$ is odd and $H:=\tilde{G}_{4,X,L}$
if $\ell=2$. Clearly, in both cases $H$ is $g$-unbounded. Thanks to Remark \ref{silver},
 there exists a surjective homomorphism
$$\rho: H \twoheadrightarrow \Gamma_L \subset \Aut(D).$$
Clearly,
$$D^H=D^{\Gamma_L}=\Q.$$
In light of Remark \ref{endo} the semisimple $\Q$-algebra
$D=\End^0(X)$ satisfies all the conditions of Theorem \ref{main0}
with $g=\dim(X)$. Applying  Theorem \ref{main0}, we conclude that
$D=\Q$, i.e., $\End^0(X)=\Q$ and therefore $\End(X)=\Z$.
\end{proof}

\section{Hyperelliptic jacobians}
\label{HJ}
\begin{thm}
\label{jac} Let $K$ be a field of characteristic zero,
 $K_a$ its algebraic closure, $f(x) \in K[x]$ a polynomial of degree $n \ge 5$ and
$\RR \subset K_a$ the set of its roots. let $K(\RR)\subset K_a$ be
the splitting field of $f$ and $\Gal(f):=\Gal(K(\RR)/K)$ the
Galois group of $f$, viewed as a subgroup of of the group
$\Perm(\RR)$ of all permutations of $\RR$. Suppose $\Gal(f)$
contains a simple non-abelian group $\mathcal{G}$ enjoying one of
the following two
 properties:

\begin{enumerate}
\item[(i)]
$n$ is odd and $\mathcal{G}$ acts $2$-transitively on $\RR$. In
addition, every $2$-minimal cover of $\mathcal{G}$ is
$\frac{n-1}{2}$-bounded.
\item[(ii)]
$n$ is even and $\mathcal{G}$ acts $3$-transitively on $\RR$. In
addition, every $2$-minimal cover of $\mathcal{G}$ is
$\frac{n-2}{2}$-bounded.
\item[(iii)]
$n$ is even and
$\End_{\mathcal{G}}(Q_{\RR})=\F_2$. In
addition, every $2$-minimal cover of $\mathcal{G}$ is
$\frac{n-2}{2}$-bounded.
\end{enumerate}

Let $J(C_f)$ be the jacobian of the hyperelliptic curve $C=C_f:
y^2=f(x)$. Then the ring $\End(J(C_f)$ of all $K_a$-endomorphisms
of $J(C_f)$ coincides with $\Z$.
\end{thm}

\begin{proof}
 Suppose $f(x)\in K[x]$ is a polynomial of degree $n \ge 5$
without multiple roots and $X:=J(C_f)$ is the jacobian of $C=C_f:
y^2=f(x)$. It is well-known that $g:=\dim(X)=\frac{n-1}{2}$ if $n$
is odd and $g:=\dim(X)=\frac{n-2}{2}$ if $n$ is even.
  It is also
well-known (see for instance Sect. 5 of \cite{ZarhinTexel}) that
$\tilde{G}_{2,X,K}\cong \Gal(f)$. More precisely, let $K(\RR)$ be
 the splitting field of $f$ and $\Gal(f):=\Gal(K(\RR)/K)$ the Galois group of $f$,
viewed as a subgroup of of the group $\Perm(\RR)$ of all
permutations of $\RR$. We have
$$\Gal(f) \subset \Perm(\RR)$$
and the action of $\Gal(f)$ on $\RR$ is transitive if and only if
$f$ is irreducible.

Now let us consider the heart (end of \S \ref{gt}) of the
permutation representation of $\Gal(f)$ on $\RR$: the faithful
$\Gal(f)$-module $Q_{\RR}$. It is well-known (see for instance,
Th. 5.1 on p. 478 of \cite{ZarhinTexel}) that  the homomorphism
$\rho_{2,X}:\Gal(K) \to \Aut(X_2)$ factors through the canonical
surjection $\Gal(K) \twoheadrightarrow \Gal(K(\RR)/K)=\Gal(f)$ and
the $\Gal(f)$-modules $X_2$ and $Q_{\RR}$ are isomorphic. In
particular,
$$\Gal(f)=\rho_{2,X}(\Gal(K))=\tilde{G}_{2,X,K}.$$

Now assume that $f$ satisfies the conditions of Theorem \ref{jac}
and let us put
$$H:=\mathcal{G} \subset \Gal(f)\subset \Perm(\RR) \subset
\Aut(Q_{\RR})=\Aut(X_2).$$
It follows easily from Lemma \ref{Kl} that we always have
$$\End_{\mathcal{G}}(X_2)=\End_H((Q_{\RR}))=\F_2.$$
Now the assertion of  Theorem \ref{jac} follows readily from Theorem \ref{apps}.
\end{proof}

\begin{proof}[Proof of Theorem \ref{mainM}]
It is well-known that all Mathieu groups $\M_n \subset \SS_n$ are,
at least, $3$-transitive permutation groups. Now Theorem
\ref{mainM} becomes an immediate corollary of Theorem \ref{jac}
combined with Lemma \ref{primer1}(a-c).
\end{proof}

\begin{proof}[Proof of Theorem \ref{mainL}]
Recall that $q$ is a power of $2$. It is well-known that
$\L_m(q)=\PSL_m(\F_q)$ acts doubly transitively on
$\P^{m-1}(\F_q)=\RR$. It is also clear that
$n=\#(\P^{m-1}(\F_q))=\frac{q^m-1}{q-1}$ is odd. Now Theorem
\ref{mainL} becomes an immediate corollary of Theorem \ref{jac}
combined with Examples \ref{seven} and Lemma \ref{Lmq}.
\end{proof}

\begin{proof}[Proof of Theorem \ref{mainL34}]
It is an immediate corollary of Theorem \ref{jac} combined with
Lemma \ref{primer1}(d) and Lemma \ref{qodd}. (Notice that in this
case one may check that the $\Gal(f)$-module $J(C_f)_2$ is very
simple.)
\end{proof}

\begin{rem}
Combining Remark \ref{primer1}(e) and Lemma \ref{Kl} with Theorem
\ref{jac}, we obtain immediately that $\End(J(C_f)=\Z$ if $n \ge
5$ and $\Gal(f)$ contains $\A_n$. This assertion was proven by a
different method (based on the very simplicity of $J(C_f)_2$) in
\cite{ZarhinMRL}.
\end{rem}

\bigskip

\noindent {\small {Department of Mathematics, Pennsylvania State University,}}

\noindent {\small {University Park, PA 16802, USA} }

\vskip .4cm

\noindent {\small {Institute for Mathematical Problems in Biology,}}

\noindent {\small {Russian Academy of Sciences, Push\-chino, Moscow Region, 142292, RUSSIA}}

 \vskip  .4cm

\noindent {\small {\em E-mail address}: zarhin@math.psu.edu}


\begin{thebibliography}{99}

\bibitem{Atlas} J. H. Conway, R. T. Curtis, S. P. Norton, R. A. Parker, R.
A. Wilson, Atlas of finite groups. Clarendon Press, Oxford, 1985.


\bibitem{FT} W. Feit, J. Tits, {\em Projective representations of
minimum degree of group extensions}. Canad. J. Math. {\bf 30}
(1978), 1092--1102.

\bibitem{Ivanov} A. A. Ivanov, Ch. E. Praeger, {\em On finite
affine 2-Arc transitive graphs}. Europ. J. Combinatorics {\bf 14}
(1993), 421--444.

\bibitem{Katz1} N. Katz, {\em Monodromy of families of curves:
    applications of some results of Davenport-Lewis}. In:
    S\'eminaire de Th\'eorie des Nombres, Paris 1979-80
    (ed. M.-J. Bertin); Progress in Math. {\bf 12},
    pp. 171--195, Birkh\"auser, Boston-Basel-Stuttgart,
    1981.

\bibitem{Katz2} N. Katz,  {\em Affine cohomological transforms,
perversity, and monodromy}. J. Amer. Math. Soc. {\bf 6} (1993),
149--222.


\bibitem{KL} P. B. Kleidman, M. W. Liebeck, {\em On a theorem of
Feit and Tits}. Proc. AMS {\bf 107} (1989), 315--322.

\bibitem{Klemm} M. Klemm, {\em \"Uber die Reduktion von
Permutationsmoduln}. Math. Z. {\bf 143} (1975), 113--117.


\bibitem{Masser} D. Masser, {\em Specialization of some hyperelliptic jacobians}. In:
 Number Theory in Progress (eds.  K. Gy\"ory, H. Iwaniec,
J.Urbanowicz), vol. I, pp. 293--307; de Gruyter, Berlin-New York,
1999.

\bibitem{Mori1} Sh. Mori, {\em The endomorphism rings of some abelian varieties}.
 Japanese J. Math,  {\bf 2}(1976), 109--130.

\bibitem{Mori2} Sh. Mori, {\em The endomorphism rings of some abelian varieties}. II.
 Japanese J. Math,  {\bf 3}(1977), 105--109.

\bibitem{Mortimer} B. Mortimer, {\em The modular permutation
representations of the known doubly transitive groups}. Proc.
London Math. Soc. (3) {\bf 41} (1980), 1--20.

\bibitem{MumfordAV} D. Mumford, Abelian varieties, Second edition,
 Oxford University Press, London, 1974.

\bibitem{Oort} F. Oort, {\em  Endomorphism algebras of abelian varieties}.
 Alebraic Geometry and Commutative Algebra  in Honor of M. Nagata (1987, Ed. H. Hijikata et al),
 Kinokuniya Cy, Tokyo 1988; Vol. II, pp. 469 - 502.

\bibitem{OZ} F. Oort, Yu. G. Zarhin, {\em   Endomorphism algebras of complex tori}.
Math. Ann. {\bf 303} (1995), 11-29.


\bibitem{Silverberg} A. Silverberg, {\em Fields of definition for homomorphisms of abelian
varieties}. J. Pure Appl. Algebra {\bf 77} (1992), 253--262.

\bibitem{SZ} A. Silverberg, Yu. G. Zarhin, {\em Variations on a theme of Minkowski and Serre}.
J. Pure and Applied Algebra {\bf 111} (1996), 285--302.

\bibitem{Suzuki} M. Suzuki, Group Theory I. Springer-Verlag, 1982.

\bibitem{TZ} Pham Huu Tiep, A. E. Zalesskii, {\em Minimal
characters of the finite classical groups}. Comm. Algebra {\bf
24}(1996), 2093--2167.


\bibitem{Wagner} A. Wagner, {\em The faithful linear representations of
$\Sn$ and $\An$ over a field of characteristic} $2$. Math. Z. {\bf
151} (1976), 127--137.



\bibitem{ZarhinMRL} Yu. G. Zarhin, {\em Hyperelliptic jacobians without
complex multiplication}. Math. Res. Letters {\bf 7} (2000), 123--132.

\bibitem{ZarhinTexel}Yu. G. Zarhin {\em Hyperelliptic jacobians and modular
representations}. In: Moduli of abelian varieties (C. Faber, G.
van der Geer, F. Oort, eds.), pp. 473--490, Progress in Math.,
Vol. 195, Birkh\"auser, Basel--Boston--Berlin, 2001.

\bibitem{ZarhinMRL2}  Yu. G. Zarhin, {\em Hyperelliptic jacobians without
complex multiplication in positive characteristic}. Math. Res.
Letters {\bf 8} (2001), 429--435.

\bibitem{ZarhinCrelle} Yu. G. Zarhin, {\em Cyclic covers of the
projective line, their jacobians and endomorphisms},
http://xxx.lanl.gov/abs/math.AG/0003002, to appear in  J. reine
angew. Math.

\bibitem{ZarhinPAMS} Yu. G. Zarhin, {\em Hyperelliptic jacobians
and simple groups} $\U_3(2^m)$. Proc. AMS, to appear.

\bibitem{ZarhinMMJ} Yu. G. Zarhin, {\em Very simple $2$-adic representations and hyperelliptic
jacobians}, http://arXiv.org/abs/math.AG/0109014.
\end{thebibliography}
\end{document}